\documentclass{article}
\usepackage{amsmath}
\usepackage{mathrsfs}
\usepackage{txfonts}
\usepackage{graphicx}
\usepackage{bm}
\usepackage{latexsym,amsfonts,amssymb}
\usepackage{enumitem}
\usepackage{multirow,hyperref,enumerate}
\setenumerate[1]{itemsep=0pt,partopsep=0pt,parsep=\parskip,topsep=5pt}
\setitemize[1]{itemsep=0pt,partopsep=0pt,parsep=\parskip,topsep=5pt}
\setdescription{itemsep=0pt,partopsep=0pt,parsep=\parskip,topsep=5pt}
\begin{document}
\setcounter{page}{1}
\newtheorem{thm}{Theorem}[section]
\newtheorem{fthm}[thm]{Fundamental Theorem}
\newtheorem{dfn}[thm]{Definition}
\newtheorem{rem}[thm]{Remark}
\newtheorem{lem}[thm]{Lemma}
\newtheorem{cor}[thm]{Corollary}
\newtheorem{exa}[thm]{Example}
\newtheorem{pro}[thm]{Proposition}
\newtheorem{prob}[thm]{Problem}
\newtheorem{fact}[section]{Fact}
\newtheorem{con}[thm]{Conjecture}
\renewcommand{\thefootnote}{}
\newcommand{\qed}{\hfill\Box\medskip}
\newcommand{\proof}{{\it Proof.\quad}}
\newtheorem{ob}[thm]{Observation}
\newcommand{\rmnum}[1]{\romannumeral #1}
\renewcommand{\abovewithdelims}[2]{%
\genfrac{[}{]}{0pt}{}{#1}{#2}}

\renewcommand{\thefootnote}{\fnsymbol{footnote}}

\title{\bf  On even-cycle-free subgraphs of the doubled Johnson graphs}

\author{Mengyu Cao\quad
Benjian Lv\footnote{Corresponding author.\ \ \newline\ \  E-mail address: caomengyu@mail.bnu.edu.cn(M.Cao), bjlv@bnu.edu.cn(B.Lv), wangks@bnu.edu.cn(K.Wang)}
\quad
Kaishun Wang\\
{\footnotesize \em  Sch. Math. Sci. {\rm \&} Lab. Math. Com. Sys.,
Beijing Normal University, Beijing, 100875,  China} }
 \date{}
 \maketitle

\begin{abstract}
The generalized Tur\'{a}n number ${\rm ex}(G,H)$ is the maximum number of edges in an $H$-free subgraph of a graph $G.$ It is an important extension of the classical  Tur\'{a}n number ${\rm ex}(n,H)$,  which is the maximum number of edges in a graph with $n$ vertices that does not contain $H$ as a subgraph. In this paper, we consider the maximum number of edges in an even-cycle-free subgraph of the  doubled Johnson graphs $J(n;k,k+1)$,  which are bipartite subgraphs of hypercube graphs. We give an upper bound for ${\rm ex}(J(n;k,k+1),C_{2r})$ with any fixed $k\in\mathbb{Z}^+$ and any $n\in\mathbb{Z}^+$ with $n\geq 2k+1.$ We also give an upper bound for ${\rm ex}(J(2k+1;k,k+1),C_{2r})$ with any $k\in\mathbb{Z}^+,$ where $J(2k+1;k,k+1)$ is known as doubled Odd graph $\widetilde{O}_{k+1}.$ This bound induces that the number of edges in any $C_{2r}$-free subgraph of $\widetilde{O}_{k+1}$ is $o(e(\widetilde{O}_{k+1}))$ for $r\geq 6,$ which also implies a Ramsey-type result.
 \vspace{2mm}


\



\end{abstract}

\section{Introduction}
Throughout this paper, all graphs  are finite undirected graphs without loops or multiple edges. Let $G=(V(G),E(G))$ be a graph with vertex set $V(G)$ and edge set $E(G).$ We use $v(G)$ and $e(G)$ to denote the number of vertices and the number of edges in $G$, respectively. For any two distinct vertices $x,y\in V(G)$, a \emph{path} of length $r$ from $x$ to $y$ in $G$ is a finite sequence of $r+1$ distinct vertices $(x=w_{0},w_{1},\ldots,w_{r}=y)$
such that $\{w_{t-1}, w_{t}\}\in E(G)$ for $t=1,2,\ldots,r$. If there is a path between any two vertices of a graph $G$, then $G$ is \emph{connected}. A cycle is a connected graph where any vertex in the graph has exactly two neighbours. A \emph{cycle} is called to be an \emph{$l$-cycle} or a \emph{cycle of length $l$} if the number of edges in the cycle is $l$, denoted by $C_l$. The phrase ``a cycle in a graph $G$'' refers to a subgraph  of $G$ which is a cycle.
Two graphs $G$ and $G'$ are \emph{isomorphic} if there is a bijection $\sigma$ from $V(G)$ to $V(G')$ such that $\{x,y\}\in E(G)$ if and only if $\{\sigma(x),\sigma(y)\}\in E(G')$.

 Let $G$ and $H$ be two graphs. We call that $G$ is \emph{$H$-free} if there does not exist a subgraph of $G$ which is isomorphic to $H$. The \emph{generalized Tur\'{a}n number} ${\rm ex}(G,H)$ is the maximum number of edges in a $H$-free subgraph of $G$. When $G=K_n$ is the complete graph of $n$ vertices, ${\rm ex}(G,H)$ is usually denoted by ${\rm ex}(n,H)$, specifying the maximum possible number of edges in an $H$-free graph on $n$ vertices. There are a huge amount of literatures investigating this function, starting with the theorems of Mantel \cite{R8} and Tur\'{a}n \cite{R9} that determine it for $H=K_r$. It is showed in \cite{R4} that  ${\rm ex}(n,H)$ is related to the chromatic number of $H$. But when $H$ is bipartite one can only deduce that ${\rm ex}(n,H)=o(n^2)$. In general, it is also a major open problem to determine the generalized Tur\'{a}n number ${\rm ex}(G,H)$ when $H$ is a bipartite graph, especially for even cycles. In this aspect, there are two widely studied functions ${\rm ex}(K_{m,n},K_{s,t})$ and ${\rm ex}(Q_n, C_{2l})$, where $K_{m,n}$ is a complete bipartite graph and $Q_n$ is a hypercube graph.

The former function ${\rm ex}(K_{m,n},K_{s,t})$, known as  the \emph{problem of Zarankiewicz} raised in 1951  (\cite{R11}),  is the analogue of Tur\'{a}n's problem in bipartite graphs. We refer the reader to \cite{Furedi-simonvits} for the details about this problem.  The latter function ${\rm ex}(Q_n, C_{2l})$, started with a problem raised by Erd\H{o}s, which is ``How many edges can a subgraph of $Q_n$ have that contains no 4-cycles?'' In \cite{R5}, Erd\H{o}s conjectured that the upper bound would be $(\frac{1}{2}+o(1))e(Q_n)$, and also asked whether $o(e(Q_n))$ edges of $Q_n$ would ensure the existence of a
cycle $C_{2l}$ for $l\geq 3.$ The best upper bound for ${\rm ex}(Q_n, C_{4})$ is obtained by Balogn et al. (\cite{Balogh2014}),  which is $(0.6068+o(1))e(Q_n)$, slightly improving the upper bounds given by Chung (\cite{R2}) and Thomason Wagner (\cite{R10}). The problem of deciding
the values of $C_6$ and $C_{10}$ is still open. In \cite{R2}, Chung showed that $ \frac{1}{4}e(Q_n)\leq {\rm ex}(Q_n, C_{6})\leq (\sqrt{2}-1+o(1))e(Q_n)$, and negatively answered the question of Erd\H{o}s for $C_6$.  Conder (\cite{Conder}) found a $3$-colouring with the same property. This implies that ${\rm ex}(Q_n, C_6) \geq\frac{1}{3}e(Q_n)$. The best upper bound is given by Balogn et al. (\cite{Balogh2014}). For some progress about ${\rm ex}(Q_n, C_{10})$, we refer the reader to \cite{Alon2007,Alon2006}. For $l\geq 2$, the upper bounds for ${\rm ex}(Q_n, C_{4l})$ and ${\rm ex}(Q_n,C_{4l+6})$ were obtained by Chung (\cite{R2}) and F\"{u}redi and \"{O}zkahya (\cite{R1}), respectively, which imply that ${\rm ex}(Q_n,C_{2l'})=o(e(Q_n))$ for $l'\geq6$ or $l'=4$. In \cite{R7}, Conlon  unified these results by showing ${\rm ex}(Q_n, H) = o(e(Q_n))$ for all $H$ that admit
a $k$-partite representation, which holds for each $H = C_{2l}$ except  $l\in \{2, 3, 5\}.$

Now we consider another noteworthy family of bipartite graphs, which are called doubled Johnson graphs. Let $n$ and $k$ be two positive integers with $n\geq k+1$. Let $[n]=\{1,2,\ldots,n\}$ and ${[n]\choose k}$ be the set of all $k$-subsets of $[n].$ The \emph{doubled Johnson graph} $J(n;k,k+1)$ is a bipartite graph with vertex set ${[n]\choose k}\cup{[n]\choose k+1},$ where two distinct vertices $u$ and $v$ are adjacent if and only if $u\subset v$ or $v\subset u.$  Recall that doubled Johnson graphs with $n=2k+1$ are usually called \emph{doubled Odd graphs}, which are distance-transitive graphs (\cite{Brouwera2}). We usually use $\widetilde{O}_{k+1}$ to denote the doubled Odd graph $J(2k+1;k,k+1)$. Notice that $J(n;k,k+1)$ is a subgraph of the hypercube $Q_n$, and the halved graphs of $J(n;k,k+1)$ are the Johnson graphs $J(n,k)$ and $J(n,k+1)$. By the definition, in the graph $J(n;k,k+1)$, the degree of each vertex in ${[n]\choose k}$ is $n-k$ and  the degree of each vertex in ${[n]\choose k+1}$ is $k+1$. Therefore, $e(J(n;k,k+1))=(n-k){n\choose k}=(k+1){n\choose k+1}.$ Since the graphs $J(n;k,k+1)$ and $J(n;n-k-1,n-k)$ are isomorphic, in the following, we only consider the case when $n\geq 2k+1$.

In this paper, we study the generalized Tur\'{a}n number ${\rm ex}(J(n;k,k+1),C_{2l})$. For each vertex $x_2$ in ${[n]\choose k+1},$ choose an edge which is incident with $x_2$. Let $E$ be the set of those edges and $K$ be the graph with vertex set ${[n]\choose k}\cup{[n]\choose k+1}$ and edge set $E$. Notice that the degree of each vertex from ${[n]\choose k+1}$ in $K$ is $1$, which implies that $K$ is cycle-free. Hence we have ${\rm ex}(J(n;k,k+1),C_{2l})\geq {n\choose k+1}=\frac{1}{k+1}e(J(n;k,k+1)).$ In the following, we consider the upper bound of ${\rm ex}(J(n;k,k+1),C_{2l})$ and obtain the following theorems.

\begin{thm}\label{thm2}
Let $k$ and $l$ be any fixed positive integers. For any $n\in\mathbb{Z}^{\rm +}$ with $n\geq 2k+1$, the following hold.
\begin{itemize}
\item[\rm(i)] For $l\geq 2$, there exists constant $c_l$ such that
\begin{align*}
{\rm ex}(J(n;k,k+1),C_{4l})\leq \left(c_l(n-k)^{-\frac{1}{2}+\frac{1}{2l}}+\frac{1}{\sqrt{k+1}}\right)e(J(n;k,k+1)).
\end{align*}
\item[\rm(ii)] For $l\geq 1$, we have
\begin{align*}
{\rm ex}(J(n;k,k+1),C_{4l+2})\leq \left(\frac{1}{2(k+1)}+\frac{\sqrt{2}}{2}+o(1)\right)e(J(n;k,k+1)),
\end{align*}
where $o(1)$ is a function $f_k(n)$ of the variable $n$ such that $\lim\limits_{n\rightarrow +\infty}f_k(n)=0.$
\end{itemize}
\end{thm}

\begin{thm}\label{thm3}
Let $l$ be a any fixed positive integer. For any $k\in\mathbb{Z}^{\rm +}$, the following hold.
\begin{itemize}
\item[\rm(i)] For $l\geq 3$, we have ${\rm ex}(\widetilde{O}_{k+1}, C_{4l})=O(k^{-\frac{1}{2}+\frac{1}{l}})e(\widetilde{O}_{k+1})$.
\item[\rm(ii)] For $l\geq 3$, we have $${\rm ex}(\widetilde{O}_{k+1}, C_{4l+2})=\left\{\begin{array}{ll}O(k^{-\frac{1}{2l+1}})e(\widetilde{O}_{k+1}),&\mbox{if}\ l=3,5,7,9,\\ O(k^{-\frac{1}{16}+\frac{1}{8(l-1)}})e(\widetilde{O}_{k+1}),&\mbox{otherwise}.\end{array}\right. $$
\item[\rm(iii)] ${\rm ex}(\widetilde{O}_{k+1}, C_{6})\leq\frac{5}{6}e(\widetilde{O}_{k+1});$ ${\rm ex}(\widetilde{O}_{k+1}, C_{8})\leq(\frac{2}{3}+o(1))e(\widetilde{O}_{k+1});$ ${\rm ex}(\widetilde{O}_{k+1}, C_{10})\leq(\frac{2}{3}+o(1))e(\widetilde{O}_{k+1}).$
\end{itemize}
\end{thm}

From Theorem \ref{thm3}, we have ${\rm ex}(\widetilde{O}_{k+1}, C_{2l})=o(e(\widetilde{O}_{k+1}))$ for $l\geq 6$, which leads to the following Ramsey-type result:
\begin{thm}\label{thm4}
Let $t$ and $l$ be positive integers with $l\geq6$. If $\widetilde{O}_{k+1}$ is edge-partitioned into $t$ subgraphs, then one of the subgraphs must contain the even cycle $C_{2l}$, provided that $k$ is sufficiently large $($depending only on $t$ and $l$$)$.
\end{thm}

This paper is organized as follows. In Section 2, we introduce some properties of the doubled Johnson graphs. In Section 3, we give an upper bound for ${\rm ex}(\widetilde{O}_{k+1},C_{2l })$ with $l=3,4,5$. In Section 4, we give an upper bound for the number of edges in $C_{4l}$-free subgraphs of $J(n;k,k+1)$ with $l\geq 2$. In Section 5, we give an upper bound for the number of edges in $C_{4l+2}$-free subgraphs of $J(n;k,k+1)$ with $l\geq 1$.

\section{Preliminary}
In this section, we will give some important properties of the doubled Johnson graphs. It is obvious that each cycle in $J(n;k,k+1)$ has even length since it is a bipartite graph.

Suppose $\Gamma$ is a graph. For any $x\in V(\Gamma),$ let $N_{\Gamma}(x)$ and $d_{\Gamma}(x)$ denote the set of neighbours of $x$ and the degree of $x$ in $\Gamma$, respectively. For any two vertices $x,y\in V(\Gamma)$, let $\partial_{\Gamma}(x,y)$ denote the distance between $x$ and $y$ in $\Gamma.$ Given a doubled Johnson graph $J(n;k,k+1)$, in the following, we usually use $V_1$ and $V_2$ to denote the set ${[n]\choose k}$ and ${[n]\choose k+1}$, respectively, which are two parts of this bipartite graph. Set $v_1:=|V_1|$ and $v_2:=|V_2|.$ Observe that $v_1={n\choose k}$ and $v_2={n\choose k+1},$ and $v_1=v_2={2k+1\choose k}$ if $n=2k+1$. For any two vertices $x$ and $y$ in $J(n;k,k+1)$, from \cite{Kong}, we have
$\partial_{\Gamma}(x,y)=|x|+|y|-2|x\cap y|.$

\begin{pro}\label{proposition1}
Let $U=(u_0,u_1,\ldots,u_i)$ be any path in $J(n;k,k+1)$. The following hold.
\begin{itemize}
\item[{\rm (i)}] If $i=3$, there exists a unique cycle of length $6$ containing $U$ in $J(n;k,k+1)$.
\item[{\rm (ii)}] If $i=2$ and $u_2\in V_1$, there exist $n-k-1$ cycles of length $6$ containing $U$ in $J(n;k,k+1)$.
\item[{\rm (iii)}] If $i=2$ and $u_2\in V_2$, there exist $k$ cycles of length $6$ containing $U$ in $J(n;k,k+1)$.
\item[{\rm (iv)}] If $i=1$, there exist $k(n-k-1)$ cycles of length $6$ containing $U$ in $J(n;k,k+1)$.
\end{itemize}
\end{pro}
\proof {\rm(i)}\ \ If $i=3$, then $u_0\in V_1$ or $u_3\in V_1$. Without loss of generality, suppose $u_0\in V_1,$ $u_0\cap u_2=F,$ $u_0=F\cup\{x\}$ and $u_2=F\cup\{y\}.$ Then $u_1=F\cup\{x,y\}.$ Assume that $u_3=F\cup\{y,z\},$ where $z\notin u_1.$ Let $w=(u_0,u_1,u_2,u_3,w_4,w_5)$ be any cycle of length $6$. Since $u_0\subseteq w_5\neq u_1$ and $|w_5\cap u_3|=k,$ we have $w_5=F\cup\{x,y\}$ and $w_4=u_3\cap w_5.$ Hence, $w$ is unique and (i) holds.

{\rm(ii)} and {\rm(iii)}\ \ By (i), it suffices to count the number of the paths $(u_0,u_1,u_2,w_3)$. If $u_2\in V_1,$ then $u_2\subseteq w_3\neq u_1$ and there are $n-k-1$ choices for $w_3.$ If $u_2\in V_2$, then $u_1\neq w_3\subseteq u_2$ and there are $k$ choices for $w_3$. Hence (ii) and (iii) hold.

{\rm(iv)}\ \ Without loss of generality, suppose $u_1\in V_1.$ There exist $n-k-1$ vertices $w_2$ such that $u_0,u_1,w_2$ is a path. By (ii), the desired result follows.    $\qed$
\begin{cor}\label{cor1.1} The following hold.
\begin{itemize}
\item[{\rm(i)}]The length of the shortest cycle in $J(n;k,k+1)$ is $6$.
\item[{\rm(ii)}]The number of $6$-cycles in $J(n;k,k+1)$ is $n(C_6)=\frac{1}{6}k(n-k-1)e(J(n;k,k+1))=\frac{1}{6}{n\choose k}(n-k)k(n-k-1).$
\end{itemize}
\end{cor}
\proof {\rm(i)} It suffices to prove that there does not exist a $4$-cycle in $J(n;k,k+1)$. Suppose $(v_1,v_2,v_3,v_4)$ is a $4$-cycle in $J(n;k,k+1)$ such that $v_1,v_3\in V_1$ and $v_2,v_4\in V_2$. Then $v_2=v_1\cup v_3=v_4,$ a contradiction.

{\rm(ii)} Since $e(J(n;k,k+1))={n\choose k}(n-k)$ and every edge is contained in $k(n-k-1)$ cycles of length $6$ by Proposition~\ref{proposition1}, we have $n(C_6)=\frac{1}{6}{n\choose k}(n-k)k(n-k-1)$.   $\qed$

In the following, we consider the number of $2$-paths in a spanning subgraph $G$ of $J(n;k,k+1).$ For any $2$-path $(x,w,y)$ in $G$, note that $x,y\in N_G(w).$ Hence, the number of $2$-paths in $G$ whose middle vertex is in $V_i$ is
\begin{align}\label{2path}
\sum_{w\in V_i}{d_G(w)\choose 2}=\frac{1}{2}\sum_{w\in V_i}d_G(w)^2-\frac{1}{2}e(G)
\end{align}
for $i=1,2$. Observe that the total number of $2$-paths in $J(n;k,k+1)$ is $\frac{n-1}{2}\cdot e(J(n;k,k+1)).$

By Cauchy-Schwarz inequality, for $i\in\{1,2\}$, we have
\begin{align}\label{CS}
\sum\limits_{w\in V_i}d_G(w)^2\geq \left(\sum\limits_{w\in V_i}d_G(w)\right)^2/v_i=e(G)^2/v_i,
\end{align}
which implies that $$\sum_{w\in V_i}{d_G(w)\choose 2}\geq\frac{1}{2v_i}e(G)^2-\frac{1}{2}e(G).$$

\section{ Upper bounds for $ex(\widetilde{O}_{k+1},C_{2l})$ with $l=3,4,5$}

Let $\mathscr{C}_6$ be the set of all $6$-cycles in $\widetilde{O}_{k+1}$ and  $G$ be any spanning subgraph of $\widetilde{O}_{k+1}.$ For any subgraphs $H$ and $L$ of $\widetilde{O}_{k+1}$, let $G\cap H$ be the graph with vertex set $V(G)\cap V(H)$ and edge set $E(G)\cap E(H)$, and $H\setminus E(L)$ be the graph with vertex set $V(H)$ and edge set $E(G)\setminus E(L)$. Notice that for any $6$-cycle $H\in\mathscr{C}_6$, $G\cap H$ is isomorphic to one of the graphs in Figure $1$. Let $\chi_0,\ \chi_1,\ \chi_2^{1},\ \chi_2^{2},\ \chi_3^{1},\ \chi_3^{2},\chi_3^{3},\ \chi_4^{1},\ \chi_4^{2},\ \chi_4^{3},\ \chi_5,\ \chi_6$ denote the ratio of the number of $6$-cycles $H$ satisfying that $G\cap H$ is isomorphic to the graphs $(1)-(12)$ in Figure $1$ to the total number of $6$-cycles in $\widetilde{O}_{k+1}$, respectively.

\begin{figure}[!h]
  \begin{center}
 \includegraphics[width=12cm]{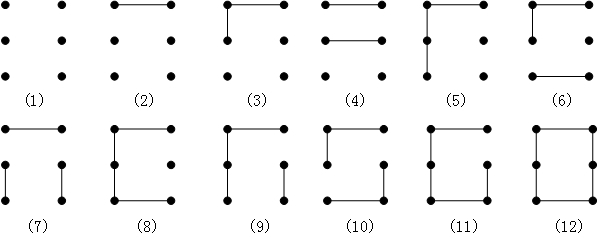}
  \end{center}
  \caption{\label{fig:InformativeFigure}Subgraphs of $C_6$}
\end{figure}
Then we have
\begin{equation}\label{chi1}
\chi_0+\chi_1+\chi_2^{1}+\chi_2^{2}+\chi_3^{1}+\chi_3^{2}+\chi_3^{3}+\chi_4^{1}+\chi_4^{2}+\chi_4^{3}+\chi_5+\chi_6=1.
\end{equation}

For any two distinct $H_1,H_2\in\mathscr{C}_6$, since the least length of a cycle in $\widetilde{O}_{k+1}$ is $6$, we have $V(H_1)\neq V(H_2)$, which implies that $G\cap H_1\neq G\cap H_2$. For any $e\in E(\widetilde{O}_{k+1}),$ let $(\mathscr{C}_6)_e$ denote the set of all $6$-cycles in $\mathscr{C}_6$ which contain $e$.  By computing the size of the set $\{(e, G\cap H)\mid H\in\mathscr{C}_6,\ e\in  E(G\cap H)\}$ in two ways, we obtain
\begin{align*}
\sum_{H\in\mathscr{C}_6}e(G\cap H)=\sum_{e\in E(G)}|(\mathscr{C}_6)_e|.
\end{align*}
By Proposition~\ref{proposition1} and Corollary~\ref{cor1.1}, we get
\begin{align}\label{chi2}
\frac{e(G)}{e(\widetilde{O}_{k+1})}&=\frac{1}{6n(C_6)}\sum_{H\in\mathscr{C}_6}e(G\cap H)\nonumber\\
&=\frac{1}{6}\left(\chi_1+2(\chi_2^{1}+\chi_2^{2})+3(\chi_3^{1}+\chi_3^{2}+\chi_3^{3})
+4(\chi_4^{1}+\chi_4^{2}+\chi_4^{3})+5\chi_5+6\chi_6\right).
\end{align}

\noindent\textbf{Proof of Theorem~\ref{thm3} (iii). } (a)\ \ Suppose $G$ is $C_6$-free. Then $\chi_6=0,$ which implies that ${\rm ex}(\widetilde{O}_{k+1},C_6)\leq\frac{5}{6}e(\widetilde{O}_{k+1})$ by (\ref{chi1}) and (\ref{chi2}).

(b)\ \ Suppose $G$ is $C_{8}$-free. For any $2$-path $L$ in $\widetilde{O}_{k+1}$, we claim that there is at most one $H$ in $\mathscr{C}_6$ such that $(G\cap H)\setminus E(L)$ is isomorphic to the graph $(8)$ in Figure $1$. Assume that $L=(u_1,u_2,u_3)$ is a $2$-path, and $H_1=(u_1,u_2,u_3,u_4,u_5,u_6)$ and $H_2=(u_1,u_2,u_3,w_4,w_5,w_6)$ are two cycles in $\mathscr{C}_6$ such that $(G\cap H_i)\setminus E(L)$ is isomorphic to the graph $(8)$ in Figure $1$ for $i\in\{1,2\}$. Since $G$ is $C_{8}$-free, $\{u_4,u_5,u_6\}\cap\{w_4,w_5,w_6\}\neq\emptyset.$ If $u_4=w_4$ or $u_6=w_6$, then $H_1$ and $H_2$ contain a same $3$-path, which is impossible by Proposition~\ref{proposition1}. If $u_4\neq w_4$ and $u_6\neq w_6$, then one can construct a cycle with length less than $6$ from $H_1$ and $H_2$. That is impossible because the least length of cycles in $\widetilde{O}_{k+1}$ is $6$. Hence, the claim holds. It is easy to see that if there is  an $H$ in $\mathscr{C}_6$ such that $(G\cap H)\setminus E(L)$ is isomorphic to the graph $(8)$ in Figure $1$, $G\cap H$ must be isomorphic to one of the graphs $(8), (11)$ and $(12)$ in Figure $1$. Conversely, for any $H\in\mathscr{C}_6$, if $G\cap H$ is isomorphic to the graph $(8), (11)$ or $(12)$ in Figure $1$, the number of $2$-paths $L$ in $\widetilde{O}_{k+1}$ such that $(G\cap H)\setminus E(L)$ is isomorphic to the graph $(8)$ in Figure $1$ is $1$, $2$ or $6$, respectively.

By counting in two ways the pairs $(L,H)$ where $L$ is a $2$-path in $\widetilde{O}_{k+1}$, $H\in\mathscr{C}_6$ such that $(G\cap H)\setminus E(L)$ is isomorphic to the graph $(8)$ in Figure $1$, we have
$$
k\cdot e(\widetilde{O}_{k+1})\geq (6\chi_6+2\chi_5+\chi_4^1)n(C_6),
$$
which implies that $\chi_4^1\leq\frac{6}{k}=o(1)$, $\chi_5\leq\frac{3}{k}=o(1)$ and $\chi_6\leq\frac{1}{k}=o(1).$ By (\ref{chi1}) and (\ref{chi2}), we have
$$
\frac{e(G)}{e(\widetilde{O}_{k+1})}\leq\frac{1}{6}(4+o(1)),
$$
and hence ${\rm ex}(\widetilde{O}_{k+1},C_{8})\leq(\frac{2}{3}+o(1))e(\widetilde{O}_{k+1})$.

(c)\ \ Suppose $G$ is $C_{10}$-free. For any $e\in E(\widetilde{O}_{k+1})$, we claim that there are at most $2k-1$ cycles $H$ in $\mathscr{C}_6$ such that $(G\cap H)\setminus \{e\}$ is isomorphic to the graph $(11)$ in Figure $1$. Assume that $e=(u_1,u_2)$ is an edge, and $H_1=(u_1,u_2,u_3,u_4,u_5,u_6)$  is a cycle in $\mathscr{C}_6$ such that $(G\cap H_1)\setminus \{e\}$ is isomorphic to the graph $(11)$ in Figure $1$. Let $H_2=(u_1,u_2,w_3,w_4,w_5,w_6)$ be any other cycle in $\mathscr{C}_6$ such that $(G\cap H_2)\setminus \{e\}$ is isomorphic to the graph $(11)$. since $G$ is $C_{10}$-free, $\{u_3,u_4,u_5,u_6\}\cap\{w_3,w_4,w_5,w_6\}\neq\emptyset.$ If $u_3\neq w_3$ and $u_6\neq w_6$, then one can construct a cycle with length less than $6$ from $H_1$ and $H_2$. That is impossible because the least length of cycles in $\widetilde{O}_{k+1}$ is $6$. Then we have $u_4=w_4$ or $u_6=w_6$. By Proposition~\ref{proposition1}, note that there are at most $2k-1$ cycles in $\mathscr{C}_6$ containing the $2$-path $(u_1,u_2,u_3)$ or $(u_6,u_1,u_2)$.  Hence, the claim holds. It is easy to see that if there is  an $H$ in $\mathscr{C}_6$ such that $(G\cap H)\setminus \{e\}$ is isomorphic to the graph $(11)$ in Figure $1$, $G\cap H$ must be isomorphic to one of the graphs $(11)$ and $(12)$ in Figure $1$. Conversely, for any $H\in\mathscr{C}_6$, if $G\cap H$ is isomorphic to the graph $(11)$ or $(12)$ in Figure $1$, the number of $e$ in $\widetilde{O}_{k+1}$ such that $(G\cap H)\setminus \{e\}$ is isomorphic to the graph $(11)$ in Figure $1$ is $1$ or $6$, respectively.

By counting in two ways the pairs $(e,H)$ where $e\in E(\widetilde{O}_{k+1})$, $H\in\mathscr{C}_6$ such that $(G\cap H)\setminus \{e\} $ is isomorphic to the graph $(11)$ in Figure $1$, we have
$$
(2k-1)\cdot e(\widetilde{O}_{k+1})\geq (6\chi_6+\chi_5)n(C_6),
$$
which implies that $\chi_5\leq\frac{6(2k-1)}{k^2}=o(1)$ and $\chi_6\leq\frac{2k-1}{k^2}=o(1).$ By (\ref{chi1}) and (\ref{chi2}), we have
$$
\frac{e(G)}{e(\widetilde{O}_{k+1})}\leq\frac{1}{6}(4+o(1)),
$$
and hence ${\rm ex}(\widetilde{O}_{k+1},C_{10})\leq(\frac{2}{3}+o(1))e(\widetilde{O}_{k+1})$. $\qed$

\section{$C_{4l}$-free subgraphs of $J(n;k,k+1)$ with $l\geq2$ }\label{even}

Let $l$ be an integer with $l\geq 2$. Suppose $G$ is a maximal spanning  $C_{4l}$-free subgraph of $J:=J(n;k,k+1)$. Notice that $V(G)=V(J)$ and $d_G(x)\geq1$ for any $x\in V(G).$

Firstly, we define an auxiliary graph $H_x:=H_x(G)$ for each vertex $x\in V(J).$ We note that the $H_x$ in this form is similar to but different from the auxiliary graph which was used by Chung \cite{R2} and F\"{u}redi et al. \cite{R1}. The vertex set of $H_x$ consists of the vertices which have distance $2$ from $x$ in $J$. In $H_x$, for any two distinct vertices $y$ and $z$, they are adjacent if and only if there exists a vertex $w\notin N_J(x)$ such that $(y,w,z)$ is a $2$-path in $G$. Notice that $|V(H_x)|=k(n-k)$ if $x\in V_1$, and $|V(H_x)|=(k+1)(n-k-1)$ if $x\in V_2$.

If $\{y,z\}\in E(H_x)$, then $\partial_G(y,z)=2.$ Hence, there exits a unique $6$-cycle containing $x,y$ and $z$ in $J$, and there exists a unique vertex $w$ such that $(y,w,z)$ is a $2$-path in $G.$ Conversely, for any two distinct vertices $y,z\in V(J)$ such that $\partial_{G}(y,z)=2$, by Proposition~\ref{proposition1} ${\rm(ii)}$ and ${\rm(iii)}$, there are $n-k-1$ (resp. $k$) vertices $x$ in $V(J)$ such $\{y,z\}\in E(H_x)$ if $y,z\in V_1$ (resp. $y,z\in V_2$). Let
$$
\mathcal{F}_i=\{(x,\{y,z\})\mid x\in V_i,\ \{y,z\}\in E(H_x)\}
$$
for  $i=\{1,2\}.$ By counting in two ways the elements in $\mathcal{F}_i,$ from (\ref{2path}), observe that
\begin{align}
&\sum\limits_{x\in V_1}e(H_x)=(n-k-1)\sum\limits_{w\in V_2}{d_G(w)\choose 2}=\frac{1}{2}(n-k-1)\sum\limits_{w\in V_2}d_G(w)^2-\frac{1}{2}(n-k-1)e(G),\label{V1eHx}\\
&\sum\limits_{x\in V_2}e(H_x)=k\sum\limits_{w\in V_1}{d_G(w)\choose 2}=\frac{k}{2}\sum\limits_{w\in V_1}d_G(w)^2-\frac{k}{2}e(G).\label{V2eHx}
\end{align}

Since $G$ is $C_{4l}$-free, we have $H_x$ is $C_{2l}$-free. If not, suppose $(y_0,y_1,\ldots,y_{2l-1})$ is a cycle in $H_x$. By the definition of $H_x,$ assume that $w_0,w_1,\ldots, w_{2l-1}\notin N_J(x)$ are the vertices such that $(y_i,w_i,y_{i+1})$ is a $2$-path in $G$ for any $i\in\{0,1,\ldots,2l-1\},$ where $y_{2l}=y_0.$ We claim that $w_0,w_1,\ldots, w_{2l-1}$ are pair-wise distinct. Suppose $w_i=w_j$ with $i\neq j.$ Then $(x,u_i,y_i,w_i),$ $(x,u_{i+1},y_{i+1},w_i)$, $(x,u_j,y_j,w_j)$ are three $3$-paths in $J,$ That is impossible since the least length of a cycle in $J$ is $6$ and there exists a unique cycle of length $6$ containing $(x,u_i,y_i,w_i)$ in $J$. Hence $(y_0,w_0,y_1,w_1,\ldots,y_{2l-2},w_{2l-2},y_{2l-1},w_{2l-1})$ is a cycle with length $4l$ in $G$, which is a contradiction. Thus,  by the consequence of Bondy and Simonovits \cite{R3}, $H_x$ can have at most $c'_l(v(H_x))^{1+1/l}$ edges, where $c'_l$ is a constant. Therefore, we have
\begin{align}\label{upperbound}
\sum_{x\in V_1}e(H_x)\leq v_l c'_1(k(n-k))^{1+1/l}.
\end{align}

\noindent \textbf{Proof of Theorem~\ref{thm2}\ (i).} Firstly, we give a lower bound of $\sum_{x\in V_1}e(H_x).$ Since $d_G(w)\leq n-k$ for any $w\in V_1,$ we get
\begin{align*}
\sum\limits_{w\in V_1}d_G(w)^2\leq (n-k)\sum\limits_{w\in V_1}d_G(w)=(n-k)e(G),
\end{align*}
which implies that
\begin{align}\label{V2ehxe}
\sum\limits_{x\in V_2}e(H_x)\leq\frac{k}{2}(n-k-1)e(G)
\end{align}
from (\ref{V2eHx}).
Since $d_G(w)^2-2d_G(w)\geq-1$ for any $w\in V_2,$ by (\ref{V1eHx}) and (\ref{V2ehxe}), we have
\begin{align}
\sum_{x\in V_1}e(H_x)-\frac{1}{k}\sum_{x\in V_2}e(H_x)
&\geq\frac{1}{2}(n-k-1)\sum_{x\in V_2}d_G(w)^2-(n-k-1)e(G)\nonumber \\
&=\frac{1}{2}(n-k-1)\left(\sum_{x\in V_2}d_G(w)^2-2\sum_{x\in V_2}d_G(w)\right) \nonumber \\
&\geq-\frac{1}{2}(n-k-1)v_2. \label{V1eHxv2}
\end{align}
Therefore, by (\ref{CS}), (\ref{V2eHx}) and (\ref{V1eHxv2}),   we have
\begin{align}\label{lowerbound1}
\sum_{x\in V_1}e(H_x)&=\left(\sum_{x\in V_1}e(H_x)-\frac{1}{k}\sum_{x\in V_2}e(H_x)\right)+\frac{1}{k}\sum_{x\in V_2}e(H_x)\nonumber\\
&\geq-\frac{1}{2}(n-k-1)v_2+\frac{1}{2}\sum\limits_{w\in V_1}d_G(w)^2-\frac{1}{2}e(G)\nonumber\\
&\geq-\frac{1}{2}(n-k-1)v_2+\frac{1}{2}\frac{e(G)^2}{v_1}-\frac{1}{2}e(G)\nonumber\\
&=-\frac{1}{2}nv_2+\frac{1}{2}\frac{e(G)^2}{v_1}.
\end{align}
By (\ref{upperbound}) and (\ref{lowerbound1}), we get
$$
e(G)^2\leq 2v_1^2c'_l(k(n-k))^{1+1/l}+nv_1v_2,
$$
which implies that
\begin{align*}
\frac{e(G)^2}{e(J)^2}&\leq 2c'_lk^{1+1/l}(n-k)^{-1+1/l}+n(n-k)^{-1}(k+1)^{-1}\\
&=(2c'_lk^{1+1/l}+k(k+1)^{-1}(n-k)^{-1/l})(n-k)^{-1+1/l}+(k+1)^{-1}
\end{align*}
from $e(J)=v_1(n-k)=v_2(k+1).$ Since $\lim\limits_{n\rightarrow+\infty}(n-k)^{-1/l}=0,$ there exists constant $c_l$ such that
$$
e(G)\leq (c_l(n-k)^{-\frac{1}{2}+\frac{1}{2l}}+(k+1)^{-\frac{1}{2}})e(J).
$$
Therefore, Theorem~\ref{thm2}\ (i) holds. $\qed$

\noindent \textbf{Proof of Theorem~\ref{thm3}\ (i).} By (\ref{CS}), (\ref{V1eHx}) and (\ref{upperbound}), we have
$$
v_1c'_l(k(k+1))^{1+1/l}\geq \frac{ke(G)^2}{2v_2}-\frac{ke(G)}{2},
$$
which implies that
$$
e(G)^2\leq 2v_1v_2c'_lk^{1/l}(k+1)^{1+1/l}+v_2e(G).
$$
Since $e(\widetilde{O}_{k+1})=v_1(k+1)=v_2(k+1),$ observe that
\begin{align*}
\frac{e(G)^2}{e(\widetilde{O}_{k+1})^2}&\leq 2c'_lk^{1/l}(k+1)^{-1+1/l}+e(G)(k+1)^{-1}e(\widetilde{O}_{k+1})^{-1}\\
&\leq(2c'_l+e(G)e(\widetilde{O}_{k+1})^{-1}(k+1)^{-2/l})(k+1)^{-1+2/l}.
\end{align*}
Thus, there exists a constant $c_l$ such that
$$
e(G)\leq c_l (k+1)^{-\frac{1}{2}+\frac{1}{l}} e(\widetilde{O}_{k+1}),
$$
and Theorem~\ref{thm3}\ (i) holds. $\qed$

\section{$C_{4l+2}$-free subgraphs of $J(n;k,k+1)$ with $l\geq 1$}

We update the auxiliary graph used in Section~\ref{even}. Let $G$ be a spanning subgraph of $J(n;k,k+1)$ and $\Omega={[n]\choose k-1}$. For any $\gamma\in\Omega$, we define a new auxiliary graph $H_{\gamma}=H_{\gamma}(G)$ as follows. The vertex set of $H_{\gamma}$ consists of all the $k$-subsets of $[n]$ which contain $\gamma$. For any two vertices $x$ and $y$ in $V(H_{\gamma}$), $x$ and $y$ are adjacent if and only if there exists a $2$-path between $x$ and $y$ in $G$.

Note that $|V(H_{\gamma})|=n-k+1$ for any $\gamma\in\Omega$. For any two distinct elements $x$ and $y$ in $V_1$, if there exists a $2$-path between $x$ and $y$ in $G$, then the $2$-path is unique in $G$, and there exists a unique $\gamma\in \Omega$ such that  $\{x,y\}\in E(H_\gamma).$ Therefore, the number of edges in $\cup_{\gamma\in\Omega}E(H_\gamma)$ equals the number of $2$-paths in $G$ whose endpoints are in $V_1$, that is
\begin{align}\label{equality1}
\sum_{\gamma\in \Omega}e(H_{\gamma})=\sum_{w\in V_2}{d_G(w)\choose 2}=\frac{1}{2}\sum_{w\in V_2}d_G(w)^2-\frac{1}{2}e(G).
\end{align}
\begin{pro}\label{cycle}
If there exists an $m$-cycle in $H_{\gamma}$ for some $\gamma\in\Omega$, then there exists a $2m$-cycle in $G$.
\end{pro}
\proof Suppose $(y_0,y_1,\ldots,y_{m})$ is a cycle in $H_{\gamma}$. By the definition of $H_{\gamma},$ assume that $w_0,w_1,\ldots, w_{m}$ are the vertices such that $(y_i,w_i,y_{i+1})$ is a $2$-path in $G$ for any $i\in\{0,1,\ldots, m\},$ where $y_{m+1}=y_0.$ We claim that $w_0,w_1,\ldots, w_{m}$ are pair-wise distinct. Suppose $w_i=w_j$ with $i\neq j.$ Then $y_i\cup y_{i+1}= y_j\cup y_{j+1}$, which is impossible since $y_i,y_{i+1},y_j$ and $y_{j+1}$ are four distinct $k$-subsets of $[n]$ which contain $\gamma$. Hence $(y_0,w_0,y_1,w_1,\ldots,y_{m-2},w_{m-2},y_{m-1},w_{m-1})$ is a cycle of length $2m$ in $G$. $\qed$

\subsection{Upper bound for ${\rm ex}(J(n;k,k+1),C_{4l+2})$ with $l\geq 1$}
\textbf{Proof of Theorem~\ref{thm2} (ii).}\ \ To get an upper bound for ${\rm ex}(J(n;k,k+1),C_{4l+2})$, we will apply the Erd\H{o}s-Stone-Simonovits Theorem \cite{R13,R12}, that if $F$ is a graph with $\chi(F)=t$ and $\chi(F\setminus \{e\})<t$ for some edge $e$ of $F$, then
$$
{\rm ex}(m,F)=\left(1-\frac{1}{t-1}+o(1)\right){m\choose 2},
$$
where $\chi(F)$ is the chromatic number of the graph $F$.

Suppose $G$ is $C_{4l+2}$-free. By Proposition~\ref{cycle}, we have $H_{\gamma}$ is $C_{2l+1}$-free.
Therefore, for $l\geq1$, according to the Erd\H{o}s-Stone-Simonovits Theorem, $H_{\gamma}$ has at most $(\frac{1}{2}+o(1)){{n-k+1}\choose 2}$ edges. By (\ref{CS}) and (\ref{equality1}), we have
$$
{n\choose k-1}\left(\frac{1}{2}+o(1)\right){{n-k+1}\choose 2}\geq \frac{e(G)^2}{2v_2}-\frac{e(G)}{2},
$$
which implies that
$$
e(G)^2\leq  v_2v_1  k(n-k)\left(\frac{1}{2}+o(1)\right)+v_2 e(G).
$$
Since $e(J)=v_1(n-k)=v_2(k+1),$ we obtain
\begin{align*}
\frac{e(G)^2}{e(J)^2}\leq\frac{k}{k+1}\left(\frac{1}{2}+o(1)\right)+\frac{1}{k+1}\frac{e(G)}{e(J)},
\end{align*}
which implies that
\begin{align*}
\frac{e(G)}{e(J)}&\leq\frac{1}{2}\left(\frac{1}{k+1}+\sqrt{\frac{1}{(k+1)^2}+\frac{4k}{k+1}\left(\frac{1}{2}+o(1)\right)}\right)\\
&\leq \frac{1}{2(k+1)}+\frac{\sqrt{1+2k(k+1)}}{2(k+1)}+o(1)\\
&\leq \frac{1}{2(k+1)}+\frac{\sqrt{2}}{2}+o(1).
\end{align*}
Therefore, Theorem~\ref{thm2} (ii) holds.  $\qed$
\subsection{$C_{4l+2}$-free subgraphs of $\widetilde{O}_{k+1}$ with $l\geq 3$}

In this subsection, let $G$ be a $C_{4l+2}$-free spanning subgraph of $\widetilde{O}_{k+1}$ with $l\geq 3$. Let $a$ and $b$ be two integers such that $4a+4b=4l+4$ and $a,b\geq2$. Notice that a cycle of length $4a$ can not intersect a cycle of length $4b$ at a single edge, otherwise their union contains a cycle of length $4l+2$. For any graph $H$, define $N(G,H)$ to be the number of subgraphs of $G$ that are isomorphic to $H$.  Firstly, we provide an upper bound on $N(G,C_{4a})$. Secondly, a lower bound on $N(G,C_{4a})$ is obtained via a lower bound on the number of $C_{2a}$'s in the auxiliary graphs constructed from $G$. Last of all, we obtain an upper bound of ${\rm ex}(\widetilde{O}_{k+1},C_{4l+2})$ and slightly improve our bound in a specific situation.

 \subsubsection{An upper bound on $N(G,C_{4a})$}
 \begin{dfn}
 The direction of an edge $\{u,v\}$ in $E(J)$, denote by $d(uv)$, to be the single number in  $u\Delta v$, where $\Delta$ is symmetric difference.
 \end{dfn}
 Let $D(F):=\{d(e)\mid e\in E(F)\}$, where $F$ is any subgraph of $\widetilde{O}_{k+1}$. Notice that for any path $P=(u_1,u_2,\ldots,u_s)$, we have $u_1\Delta u_s\subseteq D(P).$
 \begin{lem}\label{lem5.2}
 For any cycle $C$ of length $2r$ in $\widetilde{O}_{k+1}$, we have $|D(C)|\leq r.$
 \end{lem}
 \proof It suffices to prove that for any $x\in D(C)$ there exist at least two edges in $C$ whose direction is $x$. Assume that there exists $x^{\prime}\in D(C)$ such that the number of edges in $C$ with direction $x^{\prime}$ is $1$. Without loss of generality, suppose that $C=(u_1,u_2,\ldots,u_{2r})$, $d(u_{2r}u_1)=x^{\prime}$ and $x^{\prime}\in u_{2r}$. Since $x^{\prime}\notin u_1$ and $x^{\prime}\notin u_i\Delta u_{i+1}$ for $i\in\{1,2,\ldots,2r-1\},$ we have $x^{\prime}\notin u_{2r},$ a contradiction. Hence, the desired result follows. $\qed$

 \begin{lem}\label{lem1}
 Let $C$ and $C'$ be cycles of length $4a$ and $4b$ of $G$, respectively. If $E(C)\cap E(C')\neq\emptyset$, then $|D(C)\cap D(C')|\geq 2$.
 \end{lem}
\proof Suppose $\{u_1,u_2\}\in E(C)\cap E(C')$. Since $G$ has no cycles of length $4a+4b-2$, there exists $u_3\notin\{u_1,u_2\}$ such that $u_3\in V(C)\cap V(C')$. Since $|u_1\Delta u_2|=1$ and $u_3\neq u_2$, we have $u_1\Delta u_2\neq u_1\Delta u_3.$ Notice that $(u_1\Delta u_2)\cup (u_1\Delta u_3)\subseteq D(C)\cap D(C')$, which implies that  $|D(C)\cap D(C')|\geq2$. $\qed$
\begin{lem}\label{upb}
We have
$$N(G,C_{4a})= O(k^{2a-2})e(G)+O(v_1k^{2a-\frac{1}{2}+\frac{1}{b}}).$$
Moreover, if $a=b$, then $N(G,C_{4a})= O(k^{2a-2})e(G)$.
\end{lem}
\proof Let $\mathscr{C}$ denote the set of cycles of length $4a$ in $G$ and $\mathscr{C}_e$ denote the set of cycles in $\mathscr{C}$ which contain the edge $e$. Note that $|\mathscr{C}|=N(G,C_{4a})$. Let $\mathscr{E}=\cup_{C\in\mathscr{C}}E(C)$ and $\mathscr{E}:=\mathscr{E}_1\cup \mathscr{E}_2$, where $\mathscr{E}_1$ is the collection of edges that are contained in a cycle of length $4b$ in $G$, and $\mathscr{E}_2:=\mathscr{E} \setminus \mathscr{E}_1$. By counting the size of $\{(H,e) \mid H\in\mathscr{C},\ e\in\mathscr{E}\ {\mbox{and}}\ e\in E(H)\}$ in two ways, we have
\begin{align}\label{4aNe1e2}
4aN(G,C_{4a})=\sum_{e_1\in \mathscr{E}_1 }|\mathscr{C}_{e_1}|+\sum_{e_2\in \mathscr{E}_2 }|\mathscr{C}_{e_2}|.
\end{align}

Since every $4a$-cycle containing a fixed edge $e$ is determined by a sequence of directions, for each $B\in\{D(C^*) \mid C^*\in \mathscr{C}_e\}$, there are at most $|B|^{4a-1}$ $4a$-cycles $C^\prime_{4a}$ such that $D(C^\prime_{4a})=B$ and $e\in E(C^\prime_{4a}).$

For each $e_1\in\mathscr{E}_1$ (if $\mathscr{E}_1\neq\emptyset$), let $C^\prime$ be a fixed $4b$-cycle with $e_1\in E(C^\prime).$ For any $4a$-cycle $C^*\in \mathscr{C}_{e_1},$ we have $d(e)\in D(C^*)$ and $|D(C^*)\cap D(C^\prime)|\geq 2$ from Lemma~\ref{lem1}. Hence, by Lemma~\ref{lem5.2}, we have
\begin{align*}
|\{D(C^*)\mid C^*\in\mathscr{C}_{e_1}\}|\leq\sum_{i=1}^{2a-1}{|D(C^\prime)|-1\choose i}\sum_{j=0}^{2a-1-i}{2k+1-|D(C^\prime)|\choose j},
\end{align*}
which implies that
\begin{align}\label{Ce1}
|\mathscr{C}_{e_1}|\leq\sum_{i=1}^{2a-1}{D(C^\prime)-1\choose i}\sum_{j=0}^{2a-1-i}{2k+1-|D(C^\prime)|\choose j}(i+1+j)^{4a-1}=O(k^{2a-2}).
\end{align}
For each $e_2\in\mathscr{E}_2$ (if $\mathscr{E}_2\neq\emptyset$),  by Lemma~\ref{lem5.2} again, we have
\begin{align*}
|\{D(C^*)\mid C^*\in\mathscr{C}_{e_2}\}|\leq\sum_{i=0}^{2a-1}{2k\choose i},
\end{align*}
which implies that
\begin{align}\label{Ce2}
|\mathscr{C}_{e_2}|\leq\sum_{i=0}^{2a-1}{2k\choose i}(i+1)^{4a-1}=O(k^{2a-1}).
\end{align}

Notice that $|\mathscr{E}_1|\leq e(G)$ and $|\mathscr{E}_2|\leq {\rm ex}(\widetilde{O}_{k+1},C_{4b})$ because the subgraph induced by $\mathscr{E}_2$ is $C_{4b}$-free. By (\ref{4aNe1e2}), (\ref{Ce1}), (\ref{Ce2}) and Theorem~\ref{thm3} (i), we obtain
\begin{align*}
N(G,C_{4a})=\frac{1}{4a}\left(\sum\limits_{e\in \mathscr{E}_1}O(k^{2a-2})+\sum\limits_{e\in \mathscr{E}_2}O(k^{2a-1})\right)
\leq O(k^{2a-2})e(G)+O(v_1k^{2a-\frac{1}{2}+\frac{1}{b}}).
\end{align*}

In particular, if $a=b$, then $|\mathscr{E}_2|=0$. Hence,
\begin{align*}
N(G,C_{4a})=\frac{1}{4a}\sum\limits_{e\in \mathscr{E}_1}O(k^{2a-2})\leq O(k^{2a-2})e(G).
\end{align*}
We complete the proof of this lemma and obtain an upper bound of $N(G,C_{4a})$. $\qed$

 \subsubsection{A lower bound on $N(G,C_{4a})$}
 In this part, we use the auxiliary graphs defined in the beginning of this section to get a lower bound of $N(G,C_{4a})$ via a lower bound on the number of $2a$-cycles in these auxiliary graphs.

By the definition of the auxiliary graph and the proof of Proposition~\ref{cycle}, we get
\begin{equation}\label{NGNC}
N(G,C_{4a})\geq\sum_{\gamma\in\Omega} N(H_{\gamma},C_{2a}).
\end{equation}
\begin{lem}\label{thm5}{\rm(Erd\H{o}s, Simonovits~ \cite{R4})}  ~Let $L$ be a bipartite graph, where there exist vertices $x$ and $y$ such that $L\setminus\{x,y\}$ is a tree. Then for a graph $H$ with $n$ vertices and $e$ edges, there exist constants $c_1,~c_2> 0$ such that if $H$ contains more than $c_1n^{\frac{3}{2}}$ edges, then
$$
  N(H,L)\geq c_2 \frac{e^{n(L)}}{n^{2e(L)-n(L)}},
$$
where $n(L)$ and $e(L)$ are the number of vertices and  edges in $L$, respectively. $\qed$
\end{lem}

\begin{lem}\label{lowb}
We have $N(G,C_{4a})\geq c v_1\frac{\bar{d}^{4a}}{k^{2a}}-O(v_1 k^{a})$, where $\bar{d}=e(G)/v_1=e(G)/v_2$.
\end{lem}
\proof We use Lemma~\ref{thm5} with $L=C_{2a}$ in the following form so that the condition on the minimum number of edges is incorporated. Since $n(L)=e(L)=2a$, we have
\begin{align*}
N(H_\gamma,C_{2a})\geq c_2\left(\frac{e(H_\gamma)^{2a}}{v(H_\gamma)^{2a}}-\frac{(c_1v(H_\gamma)^{3/2})^{2a}}{v(H_\gamma)^{2a}}\right),
\end{align*}
which implies that
\begin{align}\label{NGNC2}
  N(G,C_{4a})\geq \sum\limits_{\gamma\in
  \Omega}c_2\left(\frac{e(H_\gamma)^{2a}}{v(H_\gamma)^{2a}}-\frac{(c_1v(H_\gamma)^{3/2})^{2a}}{v(H_\gamma)^{2a}}\right)
\end{align}
by (\ref{NGNC}).
By H\"{o}lder inequality, (\ref{CS}) and (\ref{equality1}), we have
\begin{align*}
\sum_{\gamma\in\Omega}e(H_{\gamma})^{2a}&\geq\left(\sum_{\gamma\in\Omega}e(H_{\gamma})\right)^{2a}\cdot|\Omega|^{-2a+1}=\left(\sum_{w\in V_2}{d_G(w)\choose 2}\right)^{2a}\cdot|\Omega|^{-2a+1}\\
&=\left( v_2{\bar{d}\choose 2}\right)^{2a}\cdot|\Omega|^{-2a+1}=\left( \frac{k+2}{k}{\bar{d}\choose 2}\right)^{2a}\cdot\frac{kv_1}{k+2}.
\end{align*}
Since ${\bar{d}\choose 2}/\bar{d}^2\leq\frac{1}{2}$, by (\ref{NGNC2}), we get
\begin{align*}
N(G,C_{4a})\geq c_2 v_1\frac{k}{k+2}\frac{{\bar{d}\choose 2}^{2a}}{k^{2a}}-O\left(\frac{k}{k+2}v_1 (k+2)^a\right)
\geq cv_1\frac{\bar{d}^{4a}}{k^{2a}}-O(v_1 k^a).
\end{align*}
Therefore, the desired result follows. $\qed$

 \subsubsection{Proof of Theorem~\ref{thm3} (ii)}
Since $e(G)/v_1\leq k+1$, by Lemmas~\ref{upb} and \ref{lowb}, we have
\begin{align*}
\bar{d}^{4a}&\leq O(k^{3a})+\bar{d}O(k^{4a-2})+O(k^{4a-\frac{1}{2}+\frac{1}{b}}).
\end{align*}
Hence,
\begin{align*}
\bar{d}^{4a}=\max\left\{\bar{d}O(k^{4a-2}),\ O(k^{4a-\frac{1}{2}+\frac{1}{b}})\right\},
\end{align*}
which implies that $\bar{d} = \max \left\{O(k^{1-\frac{1}{4a-1}}),\ O(k^{1-\frac{1}{4a}(\frac{1}{2}-\frac{1}{b})})\right\}$.

This bound is minimized when $a=2$ and $b=l-1$ and we get $\bar{d} =O(k^{1-\frac{1}{16}+\frac{1}{8(l-1)}})$,
which implies that
\begin{equation}\label{Oin1}
e(G) =O(k^{1-\frac{1}{16}+\frac{1}{8(l-1)}})v_1=O(k^{-\frac{1}{16}+\frac{1}{8(l-1)}})e(\widetilde{O}_{k+1}),
\end{equation}
where $e(\widetilde{O}_{k+1})=(k+1)v_1.$

Finally, we consider the case $a=b=(l+1)/2$ when $l$ is odd. By Lemmas~\ref{upb} and \ref{lowb}, we have
\begin{align*}
\bar{d}^{4a}\leq O(k^{3a})+ \bar{d}O(k^{4a-2}),
\end{align*}
which implies that $\bar{d}= O(n^{1-\frac{1}{4a-1}}).$ Since $a=b=(l+1)/2$, we immediately get $\bar{d}= O(k^{1-\frac{1}{2l+1}}),$ which implies that
\begin{equation}\label{Oin2}
e(G)=O(k^{1-\frac{1}{2l+1}})v_1=O(k^{-\frac{1}{2l+1}})e(\widetilde{O}_{k+1}).
\end{equation}

By comparing  $(\ref{Oin1})$ and $(\ref{Oin2})$ when $l$ is odd, observe that $k^{-\frac{1}{2l+1}}\leq k^{-\frac{1}{16}+\frac{1}{8(l-1)}}$ if and only if $0< l < 9.8$. Since $l\geq 3$, $(\ref{Oin2})$ improves $(\ref{Oin1})$ for $l=3,5,7,9$. We compete the proof of Theorem~\ref{thm3} (ii). $\qed$

\begin{rem}
 Our proof also implies that ${\rm ex}(\widetilde{O}_{k+1}, \Theta_{4a-1,1,4b-1})$ is $o(e(\widetilde{O}_{k+1}))$ for $a,b\geq2$ and $k\geq1$, where $\Theta_{u,v,w}$ is a theta-graph consisting of three paths of lengths $u$, $v$ and $w$ having the same endpoints and distinct inner vertices. Our result also naturally implies that $C_{2l}$ is Ramsey for $l\geq 6$, i.e., there is a monochromatic copy of $C_{2l}$ in any $t$-edge-coloring of $\widetilde{O}_{k+1}$ when $k>k(t,l)$ {\rm(}Theorem \ref{thm4}{\rm)}.
\end{rem}
\section*{Acknowledgement}
This research is supported by  NSFC (11671043).

\end{document}